\documentclass[12pt]{article}
\usepackage{amsfonts,amsmath,amssymb,bbm}
\usepackage{setspace} 
\def\proof{\noindent{\bf Proof:}\hskip10pt}  
\def\sproof{\noindent{\bf Sketch of the proof:}\hskip10pt}        
\def\QED{\hfill $\Box$}

\font\tenmath=msbm10 scaled 1200
\font\sevenmath=msbm7 scaled 1200
\font\fivemath=msbm5 scaled 1200
\newfam\mathfam \textfont\mathfam=\tenmath
\scriptfont\mathfam=\sevenmath \scriptscriptfont\mathfam=\fivemath

\begin{document}
\def \\ { \cr }
\def \R{\mathbb{R}}
\def\N{\mathbb{N}}
\def\E{\mathbb{E}}
\def\P{\mathbb{P}}
\def\Z{\mathbb{Z}}
\def\D{\mathbb{D}}
\def\C{\mathbb{C}}
\def \1{1 \mkern -6mu 1} 
\def\da{^{\downarrow}}
\def \e{{\rm e}}
\def \f{{\cal F}}
\def \tp{\tilde{\cal P}}
\def \bx{\bar{\bf x}}
\newtheorem{theorem}{Theorem}
\newtheorem{definition}{Definition}
\newtheorem{proposition}{Proposition}
\newtheorem{lemma}{Lemma}
\newtheorem{corollary}{Corollary}
\centerline{\LARGE \bf   Homogenenous Multitype Fragmentations}
\vskip 2mm

\vskip 1cm
\centerline{\Large \bf Jean Bertoin}
\vskip 1cm
\noindent
\centerline{\sl Laboratoire de Probabilit\'es et Mod\`eles Al\'eatoires}
\centerline{\sl Universit\'e Pierre et Marie Curie}
\centerline{\sl and DMA, Ecole Normale Sup\'erieure}
\centerline{\sl Paris, France}
\vskip 15mm

\noindent{\bf Summary. }{\small  A homogeneous mass-fragmentation, as it has been defined in \cite{RFC}, describes the evolution of the collection of masses of fragments of
an object which breaks down into pieces as time passes. Here, we show that this model can be enriched by considering also the types of the fragments, where a type may represent, for instance, a geometrical shape, and can take finitely many values. In this setting, the dynamics of a randomly tagged fragment play a crucial role in the analysis of the fragmentation. They are determined by a Markov additive process whose distribution depends explicitly on the characteristics of the fragmentation. As applications, we 
make  explicit the connexion with multitype branching random walks, and obtain multitype analogs of the pathwise central limit theorem and large deviation estimates for the empirical distribution of fragments.}
\vskip 3mm
\noindent
 {\bf Key words.}{ \small Multitype, fragmentation, branching process, Markov additive process.} 
 \vskip 5mm
\noindent
\noindent
{\bf A.M.S. Classification.}  {\tt 60J80, 60G18}
\vskip 3mm
\noindent{\bf e-mail.} {\tt jbe@ccr.jussieu.fr}

\begin{section}{Introduction}
 In the recent years,  there has been some interest for a class stochastic processes which are meant to serve as models for the evolution of 
an object that breaks down into smaller pieces, randomly and repeatedly as time passes. We refer to the monograph \cite{RFC} and the survey \cite{Besur} 
for a detailed account and references. Several crucial hypotheses have to be made in order to deal with models that can be analized by standard probabilistic techniques. Typically, one assumes that the process enjoys the branching property, in the sense that the dynamics of a given fragment do not depend on the others.
A further important assumption which is made in  \cite{RFC}, is that each fragment is characterized by a real number which can be viewed as its size. The latter requirement does not allow us  to consider geometrical properties like the shape of a fragment, although such notions could be relevant for describing how an object breaks down.

\begin{picture}(300,120)(-10,20)

\put (20,150){\line(1,0){100}}
\put (20,50){\line(1,0){100}}
\put (20,50){\line(0,1){100}}
\put (120,50){\line(0,1){100}}

\put (170,150){\line(1,0){100}}
\put (170,50){\line(1,0){100}}
\put (270,50){\line(-1,1){100}}
\put (170,50){\line(0,1){100}}
\put (270,50){\line(0,1){100}}
\put (220,50){\line(0,1){50}}
\put (170,100){\line(1,0){50}}

\put (320,50){\line(1,0){100}}
\put (320,50){\line(0,1){100}}
\put (320,150){\line(1,-1){100}}
\put (330,150){\line(1,-1){90}}
\put (320,150){\line(1,0){100}}
\put (420,150){\line(0,-1){100}}
\put (320,100){\line(1,0){50}}
\put (370,50){\line(0,1){50}}

\put (330,140){\line(1,0){10}}
\put (340,130){\line(1,0){10}}
\put (350,120){\line(1,0){10}}
\put (360,110){\line(1,0){10}}
\put (370,100){\line(1,0){10}}
\put (380,90){\line(1,0){10}}
\put (390,80){\line(1,0){10}}
\put (400,70){\line(1,0){10}}
\put (410,60){\line(1,0){10}}

\put (330,140){\line(0,1){10}}
\put (340,130){\line(0,1){10}}
\put (350,120){\line(0,1){10}}
\put (360,110){\line(0,1){10}}
\put (370,100){\line(0,1){10}}
\put (380,90){\line(0,1){10}}
\put (390,80){\line(0,1){10}}
\put (400,70){\line(0,1){10}}
\put (410,60){\line(0,1){10}}

\put (375,105){\line(1,0){45}}
\put (375,105){\line(0,1){45}}

\end{picture}

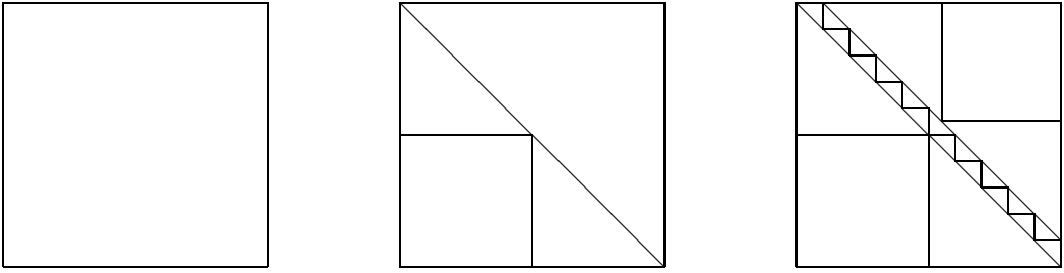
\centerline{\it Figure 1: Example of a fragmentation of a square into squares  and triangles}
\vskip 4mm

In the simpler case when time is discrete, one can analyze a fragmentation chain using the framework of branching random walks or that of multiplicative cascades. In this setting, it is therefore  natural to enrich the model by assigning to each fragment a {\it type}, which, for instance, may describe its shape, and let the evolution of each fragment depend on its initial type. The study of the latter can then be developed directly by
translating the literature on multitype branching random walks or cascades
(see e.g. \cite{Bar, BK, BRS}). 

However, we shall be interested here in the much more delicate case where time is continuous and each fragment may splits immediately, a situation which cannot be handled directly by discrete techniques based on branching processes. In the monotype setting,  Kingman's theory of exchangeable random partitions provides the key for the construction and the study of fragmentation processes in continuous time; this was pointed out first by  Pitman \cite{Pit}, see Chapter 3 in \cite{RFC} for a complete account. In the first part of this work (Sections 2 and 3), we shall briefly explain how Kingman's theory can be extended to the multitype setting (for any finite family of types), and how this extension enables us to  develop an adequate framework for multitype fragmentation processes. In short, the main result states that the dynamics of a homogeneous multitype fragmentation are characterized by a family of erosion coefficients and a family of dislocation measures. Each erosion coefficient describes the rate at which a fragment with a given type melts down as time passes, and each  dislocation measure specifies the statistics of  its sudden splits. Once the correct setting is found, statements are straightforward modifications of that in the monotype situation, and for the sake of avoiding what would be essentially a lengthy and boring duplication of existing material,
our presentation will be rather sketchy and
proofs will be omitted. Non-specialist readers may wish to consult first Chapters 2 and 3 of \cite{RFC} for getting the flavor of the arguments.

The second part of this work (Section 4) is devoted to the study of the tagged fragment, i.e. the fragment which contains a point which has been tagged at random and independently of the fragmentation process, and its applications. It departs more significantly from the monotype situation, in the sense that the evolution of the tagged fragment is now given in terms of a Markov additive process (instead of a subordinator), which depends explicitly on the characteristics of the fragmentation. The central limit theorem for Markov additive processes then enables us to determine the asymptotic behavior of certain multitype fragmentation processes, extending an old result of Kolmogorov \cite{Ko} in this area.
We will also develop the natural connexion with multitype branching random walks
from which we derive some sharp large-deviation estimates based on the work of Biggins and Rahimzadeh Sani \cite{BRS}.

Throughout this text, we shall consider a finite family of types, say with cardinal $k+1\geq 2$, which can thus be identified with $\{0,1,\ldots, k\}$. The type $0$ is special  and will only be used in peculiar situations. As it has been mentioned above, it may be convenient to think of a type as a geometrical shape (see  Figure 1 above for an example), but the type can also be used, for instance, to distinguish between active and inactive fragments in a frozen fragmentation (see \cite{DGP} for a closely related notion in the setting of coalescents). Last but not least, it was observed recently by Haas {\it et al.} \cite{HMPW} that homogeneous fragmentions bear close connexions with certain  continuum random trees, a class of random fractal spaces which has been introduced by Aldous. It is likely that more generally, multitype fragmentations can be used to construct some multifractal continuum random trees, following the analysis developed in \cite{HMPW}.

\end{section}

\begin{section}{Kingman's theory for partitions with types}

The purpose of this section is to provide a brief presentation of an extension
of Kingman's theory (see \cite{Ki} or Section 2.3.2 in \cite{RFC}) to partitions with types,
which is a key step in the analysis of random fragmentations. 

\subsection{Partitions with types}
We shall deal with two natural notions of partitions with types, which correspond to
two different points of view. The first one focuses on the masses (and the types) of the components, whereas, roughly speaking, the second one corresponds to a discretization of the object which breaks down.

We call any  numerical sequence
${\bf x}=(x_1,x_2,\ldots)$ with 
$$x_1\geq x_2\geq \cdots \geq 0\ \hbox{ and }\ \sum_1^\infty x_n\leq 1$$ 
a  {\it mass-partition}, and write ${\cal P}_{\rm m}$ for the space of mass-partitions.
A mass-partition with types is a pair $\bx=({\bf x},{\bf i})$ with ${\bf x}=(x_1,x_2,\ldots)\in{\cal P}_{\rm m}$ and ${\bf i}=(i_1,i_2,\ldots)$ a sequence in $\{0,1,\ldots, k\}$, such that for every $n\in\N
=\{1,2,\ldots\}$
\begin{equation}\label{EQ1}
x_n=0\Leftrightarrow i_n=0\,,
\end{equation}
and further
$$x_n=x_{n+1}\Rightarrow i_n\geq i_{n+1}\,,$$
i.e. the sequence $(x_1,i_1), (x_2, i_2), \ldots$ is non-increasing in the lexicographic order.  We shall write indifferently 
$$\bx=({\bf x},{\bf i})=((x_1,i_1),(x_2,i_2),\ldots)$$ 
by a slight abuse of notation.

We should think of $x_n$ as the size of the $n$-th largest component 
 of some object with total mass $1$ which has been split, and of $i_n$ as its type. A component with size $0$ means that it is absent or empty, and thus has the special type $0$.
Note that a mass-partition ${\bf x}$ can be {\it improper}, in the sense that
$\sum_1^{\infty}x_n<1$. Then the mass-defect $x_0=1-\sum_1^{\infty}x_n$ is called the mass of {\it dust}, where the dust is viewed as a set of infinitesimal particles. It may be convenient to think that the special type $0$ is also assigned to these infinitesimal particles.

We write $\bar{\cal P}_{\rm m}$ for the space of mass-partitions with types and endow it with the following distance. Let $(e_1,\ldots, e_k)$ denote the canonical basis of the Euclidean space $\R^k$, and associate to any mass-partition with types  $\bx\in\bar{\cal P}_{\rm m}$ the probability measure
on the axes of the unit cube
$$\varphi_{\bx}= x_0\delta_0+ \sum_{n=1}^{\infty} x_n \delta_{x_n e_{i_n}}\,,$$
where $x_0=1-\sum_1^{\infty} x_n$ is the mass of dust.
Then we define the distance $d(\bx,\bx')$ for every $\bx,\bx'\in\bar{\cal P}_{\rm m}$ as the Prohorov distance between the probability measures $\varphi_{\bx}$ and $\varphi_{\bx'}$,
which makes $(\bar{\cal P}_{\rm m},d)$ a compact space. We stress that
the distance $d$ is strictly weaker than other perhaps simpler distances
\footnote{It may be worthy to point out that our choice for the distance 
is well-adapted to the requirement \eqref{EQ1}. Typically, denote for $n\in\N$
by $\bx^{(n,i)}$ the mass-partition with types which consists in $n$ identical fragments 
with mass $1/n$ and fixed type $i\in\{1,\ldots,k\}$. Then $\bx^{(n,i)}$ converges as $n\to\infty$ to the degenerate partition of pure dust (and type $0$). Such a natural convergence would fail if we had chosen a stronger distance on $\bar{\cal P}_{\rm m}$ like \eqref{EQ2}.
}
 on $\bar{\cal P}_{\rm m}$ such as, for instance,
\begin{equation}\label{EQ2}
\max \{|x_n-x'_n| \1_{\{i_n=i'_n\}} + (x_n+x'_n)\1_{\{i_n\neq i'_n\}}: n\in \N\}\,.
\end{equation}

Next we turn our attention to the second notion of partition.
We call any subset of $\N=\{1,2,\ldots\}$ a block.  A {\it partition} of a block  $B\subseteq\N$ is a sequence
$\pi=(\pi_1,\pi_2,\ldots)$ of pairwise disjoint blocks with $\cup \pi_n=B$,
which is ranked according to the increasing order of the least elements,
i.e. $\inf \pi_m\leq \inf \pi_n$ whenever $m\leq n$ (with the usual convention that $\inf\varnothing=\infty$). We write ${\cal P}_{B}$ for the space of partitions of $B$.

Given a partition $\pi=(\pi_1,\pi_2,\ldots)\in{\cal P}_{B}$, we can assign to each block $\pi_n$
a type $i_n\in\{0,\ldots,k\}$, with the following convention which is related to \eqref{EQ1} :
\begin{equation}\label{EQ3}
 i_n=0\ \Leftrightarrow\ \hbox{$\pi_n$ is either empty or a singleton. }
\end{equation}
 We write $\bar\pi=(\pi,{\bf i})=((\pi_1,i_1), (\pi_2,i_2),\ldots)
 $ and call $\bar\pi$ a partition with types of $B$. We denote by $\bar{\cal P}_{B}$  the space of partitions with types of some block $B$.
 
For every block $B\subseteq \N$ and every partition $\pi=(\pi_1,\ldots)$ of $\N$,
we define $\pi_{\mid B}$, the restriction of $\pi$ to $B$, as the partition of $B$ 
whose blocks are given by $\pi_n\cap B$, $n\in\N$.
If $\bar\pi=(\pi,{\bf i})$ is now a partition with types,
we assign types to the blocks of the restricted partition $\pi_{\mid B}$ as follows.
The type of $\pi_n\cap B$ coincides with the type $i_n$ of the block $\pi_n$  if $\pi_n\cap B$ is neither empty nor a singleton, and is $0$ otherwise in order to agree with \eqref{EQ3}.
We then write $\bar\pi_{\mid B}$ for the restriction to $B$ of the partition with types $\bar \pi$.

For every pair $(\bar \pi, \bar \pi')$ of partitions with types, we  define
$$d(\bar \pi, \bar \pi')=
1/\sup\{n\in\N: \bar\pi_{\mid [n]}=\bar\pi'_{\mid [n]}\}\,,$$
where $[n]=\{1,\ldots, n\}$ and  $1/\sup \N=0$.
Note that, as the type assigned to singletons is always $0$,  the identity
$\bar\pi_{\mid [1]}=\bar\pi'_{\mid [1]}$ holds in all cases and thus 
$d(\bar \pi, \bar \pi')\leq 1$. It is easily checked that $d(\bar \pi, \bar \pi')$
 defines a distance which makes $\bar{\cal P}_{\N}$ a compact set; see Lemma 2.6 in \cite{RFC} on its page 96.

Finally, we say that a block $B\subseteq\N$ possesses an asymptotic frequency
if and only if 
$$|B|=\lim_{n\to\infty} n^{-1}{\rm Card}(B\cap [n])$$
exists. If all the blocks of $\bar\pi\in\bar{\cal P}_{\N}$ possess an asymptotic frequency,
then we say that $\bar\pi$ has asymptotic frequencies, and we write
$|\bar\pi|^{\downarrow}=({\bf x},{\bf i})$ for the sequence of the asymptotic frequencies and types of the blocks of $\bar\pi$ ranked in the non-increasing lexicographic order.
Note from Fatou's lemma that $\sum_{1}^\infty |\pi_{n}|\leq 1$ is a mass-partition and thus $|\bar\pi|^{\downarrow}
\in\bar{\cal P}_{\rm m}$.

\subsection{Exchangeability and paintbox construction}
A finite permutation is a bijection $\sigma: \N\to\N$ such that $\sigma(n)=n$
when $n$ is sufficiently large. The group of finite permutations acts naturally on the space $\bar{\cal P}_{\N}$ of partitions with types. Specifically, we write $\sigma^{-1}$ for  the finite permutation obtained as the inverse $\sigma$. Given an arbitrary partition  with types of $\N$, $\bar \pi=(\pi,{\bf i})$, $\sigma^{-1}$ maps each block $\pi_n$ of $\pi$ into a block $\sigma^{-1}(\pi_n)$ of some partition denoted by 
$\sigma(\pi)$. We decide to assign the type $i_n$ of the block $\pi_n$ to the block $\sigma^{-1}(\pi_n)$. This way, we obtain a partition with types denoted by $\sigma(\bar\pi)$. 

A measure on $\bar{\cal P}_{\N}$ is called {\it exchangeable} if it is invariant under the action of finite permutations.
Following  Kingman \cite{Ki}, we can associate to every mass-partition with types $\bx=({\bf x},{\bf i})\in\bar{\cal P}_{\rm m}$ an exchangeable probability measure  on $\bar{\cal P}_{\N}$  by the paintbox construction.
Specifically, introduce a pair of random variables $(\xi,\tau)$ with values in $\Z_+\times 
\{0,1,\ldots, k\}$ whose  distribution is specified by the following : 
$$
\P((\xi,\tau)=(n,i_n))=x_n \ \hbox { for every } n\in\N \ \hbox { and }\ 
\P((\xi,\tau)=(0,0))=1-\sum_{n=1}^\infty x_n \,.$$
Then consider a sequence $(\xi_1,\tau_1), \ldots$ of i.i.d. copies of $(\xi,\tau)$ 
and define a random partition with types $\bar\pi=(\pi,{\bf i})$ by declaring that
two distinct integers $\ell,m$ are in the same block of $\pi$
if and only if $\xi_m=\xi_{\ell}\geq 1$, and then decide that the type of that block
is $\tau_m=\tau_{\ell}$. Integers $\ell$ such that $\xi_\ell=0$ form the class of singletons of $\pi$, and their type is of course $0$.  
Similarly, if some block of $\pi$ is empty, then its type is necessarily $0$ by our convention.
The distribution of 
$\bar \pi$ will be denoted by $\rho_{\bx}$ and called the paintbox based on $\bx$. 

A slight variation of this paintbox construction  can be illustrated as follows. 
Suppose for simplicity that the mass-partition with types $\bx$ can be represented
by splitting  some geometric object, for instance a rectangle with unit area, into smaller components, for instance squares, rectangles and triangles. Each component has an area and a shape which we called a type.
Imagine that we pick at random a sequence of i.i.d. uniform points $U_1,U_2,\ldots$ in the initial object. A random partition with types is obtained by declaring that two distinct indices
are in the same block of the partition if the corresponding random points belong to the same component of the object, and the type of this block is then the type of this component.
See Figure 2 below.

\vskip5mm

\begin{picture}(400,250)(0,20)

\put (30,250){\line(1,0){400}}
\put (30,50){\line(1,0){400}}
\put (430,50){\line(0,1){200}}
\put (30,50){\line(0,1){200}}
\put (430,50){\line(-1,1){200}}
\put (230,50){\line(0,1){200}}
\put (330,50){\line(0,1){100}}
\put (330,150){\line(-1,0){300}}

\put (140,160){\makebox(0,0){$\bullet\ U_{1}$}}
\put (310,90){\makebox(0,0){$\bullet\ U_{2}$}}
\put (70,120){\makebox(0,0){$\bullet\ U_{3}$}}
\put (90,210){\makebox(0,0){$\bullet\ U_{4}$}}
\put (350,220){\makebox(0,0){$\bullet\ U_{5}$}}
\put (290,130){\makebox(0,0){$\bullet\ U_{6}$}}
\put (400,110){\makebox(0,0){$\bullet\ U_{7}$}}
\put (180,70){\makebox(0,0){$\bullet\ U_{8}$}}
\put (350,110){\makebox(0,0){$\bullet\ U_{9}$}}
\put (200,230){\makebox(0,0){$\bullet\ U_{10}$}}
\put (380,60){\makebox(0,0){$\bullet\ U_{11}$}}

\end{picture}

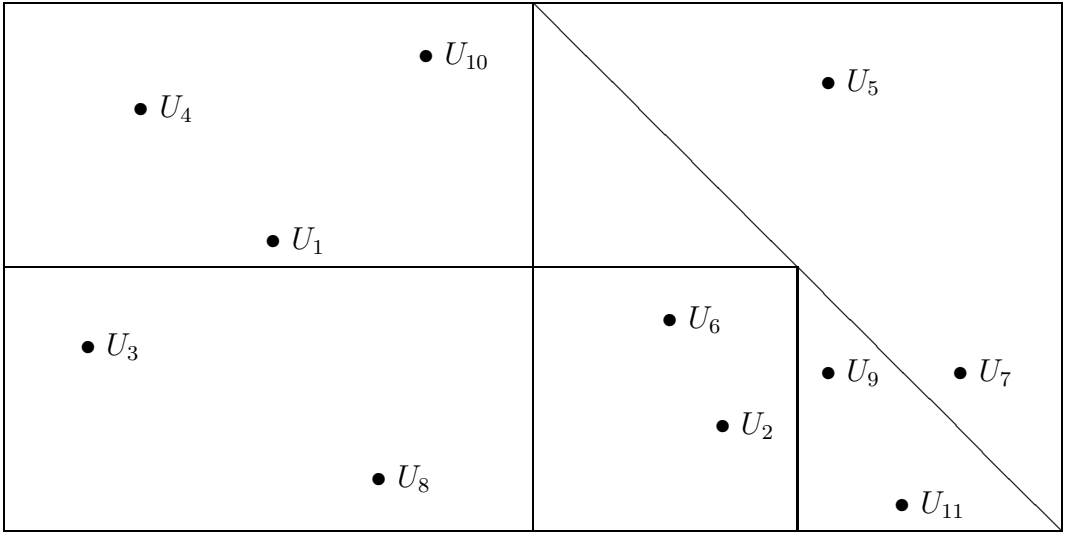
\centerline{\it Figure 2 : paintbox with types for a partition of a rectangle;}
\centerline{\it  types : square $=1$, rectangle $=2$ and  triangle $=3$}

$$\bar\pi_{\mid [11]}=((\{1,4,10\},2),(\{2,6\},1),(\{3,8\},2),(\{5,7\},3),(\{9,11\},3),(\varnothing,0),\ldots)$$

By the law of large numbers, for every positive integer $\ell$, the block $B=\{m\in\N: \xi_m= \ell\}$ has an asymptotic frequency
$$|B|=\lim_{n\to\infty} n^{-1}{\rm Card}(B\cap [n])= x_{\ell}\,,$$
and clearly the type $i_{\ell}$. One can arrange 
the sequence of the pairs (asymptotic frenquency, type) of the blocs of $\bar\pi$ in the non-increasing lexicographic order and then it coincides
with $\bx$.

Another important observation is that  $|\bar\pi_{1}|$, the asymptotic frequency and the type of the first block $\bar\pi_1$ of a paintbox based
on a mass-partition with type $\bx$, has the distribution of a size-biased sample of
$\bx$, that is
\begin{equation}\label{eqsbs}
\rho_{\bx}(|\bar\pi_1|=(x_n,i_n))= x_n\ \hbox{ for every $n\in\N$ and }\ 
\rho_{\bx}(|\bar\pi_1|=(0,0))= x_0=1-\sum_{1}^\infty x_n\,.
\end{equation}

Plainly, if $\sigma$ is a finite permutation, then $(\xi_{\sigma(1)},\tau_{\sigma(1)}),
(\xi_{\sigma(2)},\tau_{\sigma(2)}),\ldots$ is again a sequence of i.i.d. copies of $(\xi,\tau)$ and the corresponding partition with types is given by $\sigma(\bar \pi)$. Thus
$\rho_{\bx}$ is an exchangeable probability measure on $\bar{\cal P}_{\N}$, and more generally
any mixture of paintboxes produces an exchangeable probability measure on $\bar{\cal P}_{\rm m}$.
The converse to the latter assertion is a slight variation of the fundamental  theorem of Kingman \cite{Ki}, see e.g. Theorem 2.1 in \cite{RFC} on its page 100.

\begin{theorem} \label{T1}Let $\rho$ be an exchangeable probability measure on $\bar{\cal P}_{\N}$. Then $\rho$-almost every $\bar \pi\in\bar{\cal P}_{\N}$ possesses asymptotic frequencies, 
and if $\varrho$ stands for the distribution of $|\bar\pi|^{\downarrow}$  under $\rho$,
then there is the following disintegration of the measure $\rho$ : 
\begin{equation}\label{EQ4}
\rho({\rm d}\bar\pi)=\int_{\bar{\cal P}_{\rm m}}\rho_{\bx}({\rm d}\bar\pi) \varrho({\rm d}\bx)\,, \qquad \bar\pi\in\bar{\cal P}_{\N}\,.
\end{equation}

Conversely, for every probability measure $\varrho$ on $\bar{\cal P}_{\rm m}$, \eqref{EQ4}
defines an exchangeable probability measure on $\bar{\cal P}_{\N}$.
\end{theorem}

We next turn our attention to an extension of Kingman's theorem to certain sigma-finite measures on $\bar{\cal P}_{\N}$. In this direction, it is convenient to denote
for every type $i\in\{1,\ldots,k\}$ and every block $B$ that is neither empty nor a singleton, by 
${\bf 1}_{B,i}$ the partition with type of $B$ given by $((B,i), (\varnothing, 0), \ldots)$.
We also write ${\bf 1}_{i}=((1,i), (0,0), \ldots)\in\bar{\cal P}_{\rm m}$ for a related mass-partition with types.
For every $n\in\N$, we denote by $\epsilon_{(n,i)}$ for the partition with types of $\N$
which has exactly two non-empty blocks, $(\N\backslash \{n\},i)$ and $(\{n\},0)$.
The exchangeable measure on $\bar{\cal P}_{\N}$
$$\epsilon_i=\sum_{n\in\N}\delta_{\epsilon_{(n,i)}}$$
will be referred to as the erosion measure with type $i$.

The following extension of Theorem \ref{T1} to  certain possibly infinite measures is the multitype version of Theorem 3.1 in \cite{RFC} on its page 127.
Recall the notation $\bar\pi_{\mid B}$ for the partition with types restricted to some block $B$, and that $[2]=\{1,2\}$.

\begin{theorem} \label{T2}
Fix a type $i\in\{1,\ldots,k\}$ and 
let $\mu_i$ be an exchangeable measure on $\bar{\cal P}_{\N}$ such that
\begin{equation}\label{EQ5}
\mu_i(\{{\bf 1}_{\N,i}\})=0\hbox{ and }
\mu_i\left(\bar\pi\in\bar{\cal P}_{\N}: \bar\pi_{\mid [2]}\neq {\bf 1}_{[2],i}
\right)<\infty\,.
\end{equation}
Then the following
holds:

\noindent {\rm (i)} $\mu_i$-almost every partition $\bar\pi\in \bar{\cal P}_{\N}$ possesses asymptotic frequencies.

\noindent {\rm (ii)} Let
 $|\mu_i|^{\downarrow}$ be the image measure of $\mu_i$ by the mapping
$\bar\pi\to |\bar\pi|^{\downarrow}$. The restriction 
\begin{equation}\label{EQ6}
\nu_i({\rm d}{\bx})\,=\,\1_{\{{\bx}\neq{\bf 1}_i\}}|\mu_i|^{\downarrow}({\rm d}{\bx})
\end{equation}
 of
$|\mu_i|^{\downarrow}$ to $\bar{\cal P}_{\rm m}\backslash\{{\bf 1}_i\}$  fulfills 
\begin{equation}\label{EQ7}
\int_{\bar{\cal P}_{\rm m}} (1-x_1\1_{\{i_1=i\}})\nu_i({\rm d}\bx)<\infty
\end{equation}
and there is the disintegration
$$\1_{\{|\bar \pi|^{\downarrow}\neq{\bf
1}_i\}}\mu_i({\rm d} \bar\pi)\,=\,\int_{\bar{\cal P}_{\rm m}}\rho_{\bx}({\rm d}\bar\pi) \nu_i({\rm d}\bx)\,.$$

\noindent {\rm (iii)} There is a real number ${\tt c}_i\geq0$ such that
\begin{equation}\label{EQ8}
\1_{\{|\bar\pi|^{\downarrow}={\bf
1}_i\}}\mu_i({\rm d} \bar\pi)\,=\,{\tt c}_i\epsilon_i({\rm d}\bar\pi)\,.
\end{equation}

Conversely, for every real number ${\tt c}_i\geq0$ and every measure $\nu_i$  on $\bar{\cal P}_{\rm m}$
without atom at ${\bf 1}_i$  and that satisfies \eqref{EQ7}, the measure on $\bar{\cal P}_{\N}$
$$\mu_i({\rm d}\bar\pi)= {\tt c}_i\epsilon_i({\rm d}\bar\pi) + \int_{\bar{\cal P}_{\rm m}}\rho_{\bx}(d\bar\pi) \nu_i({\rm d}\bx)$$
is exchangeable and fulfills \eqref{EQ5} and \eqref{EQ6}.

\end{theorem}

As the erosion measure $\epsilon_i$ has infinite total mass, we see that ${\tt c}_i$ must be zero whenever $\mu_i$ has a finite total mass. In this situation, Theorem \ref{T2} is an immediate consequence of Theorem \ref{T1}.
 
\end{section}

\begin{section}{The structure of multitype fragmentations}
The purpose of this section is to describe the structure of multitype fragmentations.
In the monotype case, dynamics of a homogeneous fragmentation are entirely determined by an erosion rate ${\tt c}\geq 0$, which accounts for the smooth evolution of the process, and a dislocation measure $\nu$ on the space ${\cal P}_{\rm m}$ of mass-partitions, which, as its name suggests, characterizes the statistics of the sudden dislocations. 
See Sections 3.1 and 3.2 in \cite{RFC}. A similar description remains valid in the multitype situation, more precisely dynamics are then determined by a family
$({\tt c}_{i})_{i\in\{1,\ldots, k\}}$ of erosion rates and a family $(\nu_{i})_{i\in\{1,\ldots, k\}}$ of dislocation measures on $\bar{\cal P}_{\rm m}$, where the index $i$ refers of course to the type of the fragment that is eroded or dislocated.
This will be achieved first in the setting of partitions with types of $\N$, and then shifted to the more intuitive framework of mass-partitions.

\subsection{Basic definitions}

We first introduce the natural notion of homogeneous fragmentation for mass-partitions with types, which bears strong similarities with that of multitype branching process.
Specifically, let $\bar X=(\bar X(t),t\geq0)$ be a Markov
process with values in $\bar{\cal P}_{\rm m}$ and c\` adl\` ag sample paths. For every $i\in\{1,\ldots , k\}$, we write $\P_i$ for its distribution starting from $\bar X(0)={\bf 1}_{i}$, i.e. at the initial time, there is a single unit mass with type $i$.

For every mass-partition with types $\bx=({\bf x},{\bf i})$ and every real number $r\geq 0$, it will be convenient to write
$$r\bx=(r{\bf x},{\bf i})=((rx_1,i_1),(rx_2,i_2),\ldots)\,.$$
We then introduce a sequence of independent processes $Y^{(1)}, Y^{(2)},\ldots$
such that for every $n\in\N$, $Y^{(n)}$ is distributed as $x_n\bar X$ under $\P_{i_n}$.
For every $t\geq 0$, we write $Y(t)$ for the rearrangement in the non-increasing lexicographic order of the terms of the random mass-partitions with types
$Y^{(1)}(t), Y^{(2)}(t), \ldots$, and denote by $\P_{\bx}$ the distribution
of the process $Y=(Y(t), t\geq 0)$.
In particular $\P_{{\bf 1}_i}=\P_i$.

\begin{definition}\label{D2}  The process  $\bar X$ is called a
{\bf homogeneous  multitype mass-fragmentation} if, in the sense of the Markov property,  the 
law of $\bar X$ started from an arbitrary state $\bx\in\bar{\cal P}_{\rm m}$  is $\P_{\bx}$.
\end{definition}

The preceding section incites us to translate Definition \ref{D2} in the setting of partitions 
 with types of $\N$. In this direction, 
the notion of fragmentation operator (see Definition 3.1 in \cite{RFC} on its page 114)
has a natural extension in the multitype setting.

Specifically, consider $\pi\in{\cal P}_{B}$ a partition of some block  $B$ and $\bar\pi^{(\cdot)}=(\bar\pi^{(n)}, n\in\N)$ a sequence of partitions with types. We then write ${\rm Frag}(\pi,
\bar\pi^{(\cdot)})$ for the partition with types which is obtained 
from the collection of blocks with types of the sequence of the restrictions $\bar\pi^{(n)}_{\mid \pi_n}$ of $\bar\pi^{(n)}$ to the $n$-th block $\pi_n$ of $\pi$ for $n\in\N$, by rearrangement in the non-increasing lexicographic order. 
In other words, each block $\pi_n$ of $\pi$ is split using $\bar\pi^{(n)}$. 
Note that if the block $\pi_n$ is either a singleton or empty, then the partition with types $\bar\pi^{(n)}_{\mid \pi_n}$ 
does not depend on $\bar\pi^{(n)}$; more precisely, it is always given
by $((\pi_n,0),(\varnothing, 0), (\varnothing, 0), \ldots)$.

Let $\bar \Pi=(\bar\Pi(t),t\geq0)$ be a Markov
process with values in $\bar{\cal P}_{\N}$ with c\` adl\` ag sample paths. 
By a slight abuse of notation, for every $i\in\{1,\ldots , k\}$, we
write $\P_i$ for its distribution starting from $\bar\Pi(0)={\bf 1}_{\N,i}$.

\begin{definition}\label{D1}  The process  $\bar\Pi$ is called a
{\bf homogeneous multitype fragmentation} if for every time $t\geq 0$
and every type $i\in\{1,\ldots,k\}$ the distribution of $\bar \Pi(t)$ under $\P_i$
is exchangeable, and the semigroup of $\bar\Pi$  can be described as follows :
 
 Fix $t,t'\geq0$ and consider a partition with types $\bar\pi=(\pi,{\bf i})$  where ${\bf i}=(i_1,i_2,\ldots)$ is a sequence in $\{0,1,\ldots, k\}$. Let
$\bar\pi^{(\cdot)}=(\pi^{(1)}, \ldots)$ denote a sequence of independent exchangeable
random partitions with types, such that for every $n\in\N$ with $i_n\neq 0$,  $\bar\pi^{(n)}$ is distributed as $\bar\Pi(t')$ under $\P_{i_n}$. When $i_n=0$, the block $\pi_n$ 
is either empty or a singleton; the role of  $\bar\pi^{(n)}$  has no importance and its law can be chosen arbitrarily. 
 The conditional distribution of $\bar \Pi(t+t')$ given
$\bar \Pi(t)=(\pi,{\bf i})$ is then the law of ${\rm Frag}(\pi,\bar\pi^{(\cdot)})$.
\end{definition}

Let us now explain the connexion between these two definitions.
When $\bar\Pi$ is a homogeneous  multitype fragmentation, we know from Kingman's Theorem \ref{T1} that for every $t\geq 0$, the exchangeable random partition with types $\bar\Pi(t)$ possesses asymptotic frequencies a.s. If we write
$\bar X(t)=|\bar\Pi(t)|^{\downarrow}$ for the random multitype mass-partition obtained by reordering these asymptotic frequencies in the non-increasing lexicographic order,
it can be proved that the process $\bar X=(\bar X(t), t\geq 0)$ is then a homogeneous 
multitype mass-fragmentation. Technically, the main difficulty is to establish
that the paths of $t\to |\bar\Pi(t)|^{\downarrow}$ are c\`adl\`ag; in this direction we stress that this could fail if we had equipped $\bar{\cal P}_{\rm m}$ with a stronger distance such as
that given by \eqref{EQ2}.
 In the converse direction, one can rephrase the argument of Berestycki \cite{Ber1} and show that given a homogeneous multitype mass-fragmentation
$\bar X=(\bar X(t), t\geq 0)$, there exists a homogeneous multitype fragmentation 
$\bar\Pi$ such that the process $(|\bar\Pi(t)|^{\downarrow}, t\geq 0)$ is distributed as $\bar X$. In short,  there a bijective correspondence between the laws of homogeneous  multitype mass-fragmentations and laws of homogeneous multitype fragmentations.

Of course, the fundamental difference with the monotype case (see for instance Definition 3.2 in \cite{RFC} on its page 119) is that the distribution of the sequence $\bar\pi^{(\cdot)}$ which is used to split the partition with types $\bar \Pi(t)$ into finer the blocks depends on $\bar \Pi(t)$. However, this dependence only arises through the types of the blocks of $\bar \Pi(t)$, and does not 
involve directly  the partition $\Pi(t)$. This preserves the possibility of adapting 
the approach developed in Chapter 3 of \cite{RFC}, provided that one can handle some technical issues.

In particular, it is easily seen that the fragmentation operation is compatible with the restriction of partitions with types, in the sense that for every integer $n$
\begin{equation}\label{EQ9}
{\rm Frag}(\pi,\bar\pi^{(\cdot)})_{\mid [n]}\,=\,{\rm Frag}(\pi_{\mid
[n]},\bar \pi^{(\cdot)})\,.
\end{equation}
This entails that the Markov property still holds for the restricted process $\bar\Pi_{\mid [n]}= (\bar\Pi_{\mid [n]}(t),t\geq0)$, and as the latter can only take finitely many values,
$\bar\Pi_{\mid [n]}$ is a Markov chain in continuous times. Note that
$\bar\Pi_{\mid [n]}$ coincides with the restriction of $\bar\Pi_{\mid [n+1]} $ to $[n]$,
and that the initial process $\bar\Pi$ can be recovered from the sequence of Markov chains $\bar\Pi_{\mid [n]}$, $n\in\N$.

Just as in Section 3.1.2 of \cite{RFC}, these observations enable us to characterize the law of $\bar\Pi$ by a finite family of measures $(\mu_i)_{ i\in\{1,\ldots, k\}}$ on $\bar{\cal P}_{\N}$ as follows.
For every fixed $i\in\{1,\ldots, k\}$, $n\geq 2$ and every partition with types $\bar\gamma$
of $[n]$ with $\bar\gamma\neq {\bf 1}_{[n],i}$, we introduce the jump rate of $\bar\Pi_{\mid [n]}$ from ${\bf 1}_{[n],i}$ to $\bar\gamma$,
$$q_{n,\bar\gamma}\,=\,\lim_{t\to0+}{1\over t}\P_i\left(\bar\Pi_{\mid
[n]}(t)=\bar\gamma\right)\,.$$
By the very same arguments as in Section 3.1.2 of \cite{RFC}, one can check that the collection of those jump rates characterize the evolution of the restricted Markov chains $\bar\Pi_{\mid [n]}$, and thus of the process $\bar \Pi$.
Further these jump rates can be represented  as
\begin{equation}\label{EQ10}
q_{n,\bar\gamma}=\mu_i\left(\left\{\bar\pi\in\bar{\cal P}_{\N}: \bar\pi_{\mid [n]}=\bar\gamma
\right\}\right)\,,
\end{equation}
where $\mu_i$  is an exchangeable measure on $\bar{\cal P}_{\N}$ with
$$
\mu_i(\{{\bf 1}_{\N,i}\})=0 \hbox{ and } \mu_i\left(\left\{\bar\pi\in\bar{\cal P}_{\N}: \bar\pi_{\mid [n]}\neq {\bf 1}_{[n],i}
\right\}\right)<\infty \hbox{ for any }n\geq 2\,,$$
and these requirements determine the measure $\mu_i$ uniquely.
Note that when the condition above is fulfilled for $n=2$, then, thanks to the exchangeability, it is fulfilled for every $n\geq 2$, therefore it is equivalent
to \eqref{EQ5}. 
We shall refer to the family $(\mu_i)_{ i\in\{1,\ldots, k\}}$ as the splitting rates of $\bar\Pi$.

\subsection{Poissonian constructions}
Our first goal in this section is to show that  any family of exchangeable measures $(\mu_i)_{ i\in\{1,\ldots, k\}}$ 
which fulfill \eqref{EQ5} can be viewed as the splitting rates of a homogeneous multitype  fragmentation $\bar\Pi$. More precisely, we shall briefly present a Poissonian construction of $\bar\Pi$ which mimics that in Section 3.1.3 of \cite{RFC} in the monotype case.
For the sake of simplicity, we assume that the initial state has been chosen equal to ${\bf 1}_{\N,i_{0}}$ for some type $i_{0}\in\{1,\ldots,k\}$.

For every type $i\in\{1,\ldots, k\}$, consider the atoms $(t_{i,m}, \bar\pi^{(i,m)}, \ell_{i,m})_{m\in\N}$
of a Poisson random measure 
in $\R_+\times \bar{\cal P}_{\N}\times \N$
with intensity
${\rm d}t \otimes \mu_i\otimes \#$,  where $\#$ stands for the counting measure on $\N$. This means that for every measurable set $A\subseteq \R_+\times \bar{\cal P}_{\N}\times \N$,
the cardinal of the collection of indices $m$ for which $(t_{i,m}, \bar\pi^{(i,m)}, \ell_{i,m})\in A$
has the Poisson distribution with parameter ${\rm d}t \otimes \mu_i\otimes \#(A)$
and to disjoint sets correspond independent Poisson variables.
We assume that these Poisson random measures are independent for different types $i\in\{1,\ldots,k\}$.

For every integer $n$, we can construct a Markov chain in continuous time $\bar\Pi^{[n]}=(\bar\Pi^{[n]}(t), t\geq 0)$ with values in $\bar{\cal P}_{[n]}$
as follows. The atoms $(t,\bar\pi, \ell)$ of the Poisson point measures such that $\bar\pi_{\mid [n]}
= {\bf 1}_{[n],i}$ or $\ell> n$ play no role in the construction and can thus be removed.
Thanks to \eqref{EQ5}, the instants $t$ at which an atom $(t,\bar\pi,\ell)$ that has not been removed  arises, form a discrete set of $\R_+$,
and the chain $\bar\Pi^{[n]}$ can only jump at such times.
More precisely, if $(t_{i,m}, \bar\pi^{(i,m)}, \ell_{i,m})$ is an atom which has not been removed, then we look at the $\ell_{i,m}$-th block of 
$\bar\Pi^{[n]}(t_{i,m}-)$, say $B$. If the type of this block is different from the type $i$ of the atom (in particular if $B$ is either empty or a singleton), then we decide that $\bar\Pi^{[n]}(t_{i,m})=\bar\Pi^{[n]}(t_{i,m}-)$ and $t_{i,m}$ is not a jump time for the chain  $\bar\Pi^{[n]}$. If the type of $B$ is the same as the type $i$ of the atom, then
$\bar\Pi^{[n]}(t_{i,m})$ is the partition obtained
from $\bar\Pi^{[n]}(t_{i,m}-)$ by replacing $B$, that is the $\ell_{i,m}$-th block of $\bar\Pi^{[n]}(t_{i,m}-)$,  by the restriction of $\bar\pi^{(i,m)}$ to this block, and leaving the other blocks and types unchanged.

To give an example, take for instance $n=7$,  $\ell_{i,m}=2$,
$$\bar\pi^{(i,m)}=\left( (\{1,3,5,7, \ldots\},3),(\{2,4,6, \ldots\},1), \ldots\right)\,,$$
and set for simplicity $t=t_{i,m}$.
Assume also that
$$\bar\Pi^{[n]}(t-)= \left( (\{1,2\},1),(\{3,4,5\},i), (\{6,7\},2), (\varnothing, 0), \ldots\right)\,.$$
As $\ell_{i,m}=2$, we look at the 2nd block of $\bar\Pi^{[n]}(t-)$, which is $B=\{3,4,5\}$ and has type $i$, and thus coincides with the type the atom $(t_{i,m}, \bar\pi^{(i,m)}, \ell_{i,m})$. At time $t$,  we split $B$ using the partition with types $\bar\pi^{(i,m)}$.
This produces two new blocks with types:  $\{3,5\}$ which has type $3$, and $\{4\}$
which is a singleton and thus has type $0$. We conclude that
$$\bar\Pi^{[n]}(t)= \left( (\{1,2\},1),(\{3,5\},3), (\{4\},0), (\{6,7\},2), (\varnothing, 0), \ldots\right)\,.$$

It is easily seen that this construction is compatible with the restriction, in the sense that
for every $n\in\N$, 
$\bar\Pi^{[n]}$ coincides with the restriction of $\bar\Pi^{[n+1]}$ to $[n]$.
We refer to Lemma 3.3 in \cite{RFC} on its page 118 for  the argument in the monotype case. 
This implies that there exists a process $\bar\Pi$ with values in $\bar{\cal P}_{\N}$ such that 
the restriction of $\bar\Pi$ to $[n]$ is $\bar\Pi^{[n]}$; see Lemma 2.5 in \cite{RFC} on its page 95 for a closely related argument.

A crucial step  is to show that for every $t\geq 0$, the distribution of  $\bar\Pi(t)$
 is exchangeable.
The proof relies on the following technical lemma, which is the multitype version of Lemma 3.2 in \cite{RFC} on its page 116.

\begin{lemma}\label{L1}  Let $\bar\pi=(\pi,{\bf i})\in \bar{\cal P}_{\N}$ be an  exchangeable random 
partition with types and  $\bar\pi^{(\cdot)}=(\bar\pi^{(n)}, n\in \N)$ a sequence of random partitions with types. Suppose that :

\noindent $\bullet$  $\pi$ and  $\bar\pi^{(\cdot)}$ are  independent conditionally on ${\bf i}$,

\noindent $\bullet$ the sequence  $\bar\pi^{(\cdot)}$ is
{\bf doubly-exchangeable}\index{doubly-exchangeable}, in the sense that for every finite permutation $\sigma$
of $\N$, the sequences 
$$\left(\sigma(\bar \pi^{(n)}), n\in \N\right)\quad
\hbox{and}\quad 
\left(\bar\pi^{(\sigma(n))}, n\in \N\right)$$
both have the same law as 
$\bar\pi^{(\cdot)}$.
Then the random partitions with types  $\bar\pi$ and ${\rm Frag}(\pi,\bar \pi^{(\cdot)})$ are
jointly exchangeable, that is their joint distribution is
invariant by the action of permutations.
\end{lemma}

\sproof One observes that with probability one, the conditional distribution of $\pi$ given ${\bf i}$ is an exchangeable probability measure on ${\cal P}_{\N}$. One can then follow the argument of the proof of Lemma 3.2 in \cite{RFC}. \QED

It is then easy to verify from standard properties of Poisson random measures that 
the process $\bar\Pi^{[n]}$ which has just been constructed is a Markov chain in continuous time, and that its jumps rates 
$$q_{n,\bar\gamma}\,=\,\lim_{t\to0+}{1\over t}\P_i\left(\bar\Pi^{[n]}(t)=\bar\gamma\right)
$$
for every $\bar\gamma\in \bar{\cal P}_{[n]}$ with $ \bar\gamma\neq {\bf 1}_{[n],i}$, are given by \eqref{EQ10}.
This shows that the process $\bar\Pi$, which is specified by the requirement that
its restriction to $[n]$ coincides with $\bar\Pi^{[n]} $, is a homogeneous multitype fragmentation
with  splitting rates  $(\mu_i)_{ i\in\{1,\ldots, k\}}$.
Applying Theorem \ref{T2}, we can summarize this analysis in the following statement.

\begin{proposition}\label{P1} Let $\bar \Pi$ be a homogeneous multitype fragmentation.
There exists a unique family 
$({\tt c}_i)_{ i\in\{1,\ldots, k\}}$ of nonnegative real numbers and a unique family 
$(\nu_i)_{ i\in\{1,\ldots, k\}}$ of measures on $\bar{\cal P}_{\rm m}$ which fulfill \eqref{EQ7},
such that the family $(\mu_i)_{ i\in\{1,\ldots, k\}}$ of splitting rates of $\bar\Pi$
is given by \eqref{EQ8}.

Conversely, for every family 
$({\tt c}_i)_{ i\in\{1,\ldots, k\}}$ of nonnegative real numbers and every family 
\linebreak
$(\nu_i)_{ i\in\{1,\ldots, k\}}$ of measures on $\bar{\cal P}_{\rm m}$ which fulfill \eqref{EQ7},
if we define measures $\mu_i$ on $\bar{\cal P}_{\N}$  by \eqref{EQ8},
then the Poissonian construction above produces a homogeneous multitype fragmentation $\bar \Pi$ with splitting rates $(\mu_i)_{ i\in\{1,\ldots, k\}}$.

\end{proposition}

Berestycki \cite{Ber1} established a related Poissonian construction for monotype  mass-fragmentations. The latter can be extended to the multitype setting provided that the erosion coefficients ${\tt c}_i$  are all the same, which enlighten the probabilistic interpretation of the dislocation measures $\nu_i$. 

Specifically, 
for each type $i\in\{1,\ldots,k\}$, consider the atoms
$(t_{i,m}, \bx^{(i,m)}, \ell_{i,m})_{m\in\N}$
of a Poisson random measure 
in $\R_+\times \bar{\cal P}_{\rm m}\times \N$
with intensity
${\rm d}t \otimes \nu_i\otimes \#$,  where $\#$ stands for the counting measure on $\N$.
Assume as usual that these Poisson measures are independent for different types.
One can construct a pure jump process $\left(\bar Y(t), t\geq0\right)$
in $\bar{\cal P}_{\rm m}$ which
jumps only at times $t_{i,m}$ at which some atom $(t_{i,m}, \bx^{(i,m)}, \ell_{i,m})$ occurs. The jump (i.e. the dislocation) induced by such an atom
can  be described as follows. 

We consider the mass-partition with types immediately before time $t_{i,m}$, that is 
$\bar Y(t_{i,m}-)$, and look at its $\ell_{i,m}$-th term, say $(y,j)$ for some $y\geq 0$ and
$j\in\{0,\ldots, k\}$
(recall that the terms of a mass-partition with types are ranked in the non-increasing lexicographic order).
If the type $j$  is different from the type $i$ of the atom, then we simply set $\bar Y(t_{i,m}-)=\bar Y(t_{i,m})$. Otherwise, the 
$\ell_{i,m}$-th term of $\bar Y(t_{i,m}-)$ is dislocated according to $\bx^{(i,m)}$, that is it is replaced by the mass-partition with types $y\bx^{(i,m)}$. The other terms of $\bar Y(t_{i,m}-)$ are left unchanged, and  $Y(t_{i,m})$ then results from
the rearrangement in the non-increasing lexicographic order of all the terms.

For instance, if 
$$\ell_{i,m}=2\ ,\ \bx^{(i,m)}=((\frac{2}{3},2),(\frac{1}{3},1),(0,0),\ldots)$$ and
$$\bar Y(t_{i,m}-)=\left((\frac{1}{2},4),(\frac{1}{3},i),(\frac{1}{6},1),(0,0),\ldots\right)$$
then at time $t_{i,m}$ the second term of $\bar Y(t_{i,m}-)$, i.e. $(\frac{1}{3},i)$ is dislocated using $\bx^{(i,m)}$. This
produces the sequence $((\frac{2}{9},2),(\frac{1}{9},1),(0,0),\ldots)$, and finally
$$\bar Y(t_{i,m})=\left((\frac{1}{2},4),(\frac{2}{9},2),(\frac{1}{6},1),(\frac{1}{9},1),(0,0),\ldots\right)\,.$$

The process $\bar Y$ is then a homogeneous multitype fragmentation with
zero erosion and dislocation measures $(\nu_i)_{i\in\{1,\ldots,k\}}$.
Following an argument in Berestycki \cite{Ber1}, one can check that for every ${\tt c}\geq 0$, the exponentially discounted process $\bar X(t)=\e^{-{\tt c}u}\bar Y(t)$, $t\geq 0$, 
is then a homogeneous multitype mass-fragmentation with dislocation measures
$(\nu_i)_{i\in\{1,\ldots,k\}}$ and erosion coefficients ${\tt c}_i={\tt c}$ for every 
$i\in\{1,\ldots,k\}$. Unfortunately, this simple transformation cannot be extended to the case when the erosion coefficients are distinct. 
Informally, when the erosion coefficients depend on the type of the fragments, one would need information about the types of the ancestors of each fragment of $\bar Y(t)$ 
 in order to determine
the proportion of its mass that has been turned to dust at time $t$.
This information is available for processes with values in $\bar{\cal P}_{\N}$, but not for those with values in $\bar{\cal P}_{\rm m}$.

\end{section}

\begin{section}{The tagged fragment}
Up to a few technical issues, the analysis of multitype fragmentations was so far an easy translation of that in the monotype situation. However more significant differences appear when dealing with finer aspects of these processes.
Here, we shall focus on the evolution of the tagged fragment, that is the fragment which contains a point which has been picked uniformly at random and independently of the fragmentation process. The relevance of this study stems from the fact that,
even though the tagged fragment alone does not characterize the evolution of the fragmentation, it captures some useful information. In particular, this will enable us to determine the asymptotic behavior of the fragmentation, by making explicit the connexion with multitype branching random walks.
 
Let  $\bar{\bf X}=(X,T)$ be a homogeneous multitype mass-fragmentation, where $X_n(t)$ stands for the mass of the $n$-th largest fragment at time $t$ and $T_n(t)$ for its type. It will be convenient to think of $\bar{\bf X}$ as associated to 
a homogeneous multitype fragmentation $\bar\Pi$ by $\bar{\bf X}=|\bar\Pi|^{\downarrow}$.
 In order to avoid technical discussions, we shall  assume throughout this section that the fragmentation process is {\it conservative}, i.e. for every type $i\in\{1,\ldots, k\}$, the erosion coefficient is ${\tt c}_i=0$
and the dislocation measure satisfies
\begin{equation}\label{eqcons}
\nu_i\left(\left\{\bx\in\bar{\cal P}_{\rm m}: x_0=1-\sum_1^{\infty}x_n>0\right\}\right)=0\,.
\end{equation}

 The description of the evolution of the tagged fragment relies on the notion of {\it Markov additive processes}. We first provide some background in this area, referring to Section XI.2 of Asmussen \cite{As} for details.

\subsection{Background on Markov additive processes}
The class of Markov additive processes that will be useful in this work is that formed by bivariate Markov processes $(J_t,S_t)_{t\geq 0}$,
where $(J_t)_{t\geq 0}$ is a continuous time Markov chain with values in the finite space of types $\{1,\ldots,k\}$, and, roughly speaking, on every time-interval on which 
$J$ stays constant, $S$ evolves as a subordinator (i.e. an increasing process with independent and stationary increments) with characteristics specified by the value of $J$.
More precisely, one requires that for every $t,t'\geq 0$,
\begin{equation}\label{eqmap}
\E((f(S_{t+t'})-f(S_t))g(J_{t+t'})\mid  J_{t},S_{t})=\E_{J_t,0}(f(S_{t'})g(J_{t'}))
\end{equation}
where $\E_{j,s}$ refers to the mathematical expectation when the process $(J,S)$ starts from the state $(j,s)$ and $f,g$ denote two generic measurable nonnegative functions.

The law of the  Markov  chain $(J_t)_{t\geq 0}$ is specified by its intensity matrix 
${\bf \Lambda}=(\lambda_{ij})_{i,j\in\{1,\ldots,k\}}$, i.e. for $i\neq j$, $\lambda_{ij}$ is the jump rate of $J$ from $i$ to $j$ and $\sum_{j=1}^k \lambda_{ij}=0$. On every time-interval $[t,t+t'[$ on which $J\equiv i$, 
$S$ evolves as a subordinator $S^{(i)}$ with Bernstein
\footnote{
Note that, since we are dealing with subordinators,  we shall work with the Bernstein exponent $\Psi^{(i)}$ whereas Asmussen \rm \cite{As} uses the cumulant
$\theta\to -\Psi^{(i)}(-\theta)$.} 
 exponent $\Psi^{(i)}$, i.e.
$$\E(\exp(-\theta S^{(i)}_t))=\exp(-t\Psi^{(i)}(\theta))\,.$$
The Bernstein exponent is a concave increasing function which can take the value  $-\infty$, and is nonnegative and finite on $[0, \infty[$.

Further, a jump of $J$ from $i$ to $j\neq i$ has a probability $p_{ij}$ of inducing a jump of $S$ at the same time, the distribution of which is denoted by $B_{ij}$,
and we write 
$\hat B_{ij}(\theta)=\int \e^{-\theta x} B_{ij}({\rm d}x)$. It is also convenient to agree that $p_{ii}=0$.

For every types $i,j\in\{1,\ldots,k\}$ and every $\theta\geq 0$ and $t\geq 0$, there is the following identity between $k\times  k$ matrices :
$$\E_{i,0}(\e^{-\theta S_t}, J_t=j)= \left(\e^{-t {\bf \Phi}(\theta)}\right)_{{ij}} \,,$$
where
\begin{equation}\label{eqK}{\bf \Phi}(\theta) = -{\bf \Lambda}+\left(\Psi^{(i)}(\theta)\right)_{\rm diag}+(\lambda_{ij}p_{ij}
(1-\hat B_{ij}(\theta)))\,,
\end{equation}
see Proposition 2.2 in \cite{As} on its page 311. We shall refer to ${\bf \Phi}$ as the {\it Bernstein matrix} of $(J,S)$.

\subsection{Distribution of the tagged fragment}
We are interested in the process of the asymptotic frequency and the type of the first block $(|\bar\Pi_1(t)|, t\geq 0)$ of a homogeneous  multitype fragmentation $\bar\Pi$.
Recall from the paintbox construction that $|\bar\Pi_1(t)|$ can be viewed as the mass and type of the fragment which contains some point that has been picked at random according to the mass-distribution and independently of the fragmentation.

The conditions which have been enforced at the beginning of this section ensure that for every $t\geq 0$, the first block
$\Pi_1(t)$ is neither empty nor a singleton, hence its asymptotic frequency is strictly positive and its type is not $0$.
This allows us to introduce the process $(J,S)$ with values in $\{1,\ldots,k\}\times \R_+$
by
 $$|\bar\Pi_1(t)|=(\exp(-S_t), J_t)\,,Ê\qquad t\geq 0\,.$$
 
\begin{theorem}\label{T3} Suppose that the homogeneous  multitype fragmentation $\bar\Pi$ has erosion coefficients ${\tt c}_i=0$
and that its dislocation measures fulfill \eqref{eqcons}. 
Then $(J,S)$ is a Markov additive process with Bernstein matrix
given for every $\theta\geq 0$ by
$$
{\bf \Phi}(\theta) =\left(\int_{\bar{\cal P}_{\rm m}}\left({\1}_{\{i=j\}}-\sum_{n=1}^\infty x_n^{1+\theta}{\1}_{\{i_n=j\}}\right)\nu_i({\rm d}{\bx})\right)_{i,j\in\{1,\ldots,k\}}\,.
$$
\end{theorem}
\proof The fact that $(J,S)$ is a Markov process that satisfies \eqref{eqmap} can be seen from the Poissonian construction and the arguments in Section 3.2.2 of \cite{RFC}.

The determination of the Bernstein matrix also relies on the Poissonian construction.
First, note that for $i\neq j$, the jump rate $\lambda_{ij}$ of the type process $J$ coincides with the rate of occurrence of atoms $(t_{i,m},\bar\pi^{(i,m)},1)$ with $\bar\pi^{(i,m)}=(\pi^{(i,m)},{\bf i})$
and $i_1=j$. Using \eqref{eqsbs} and Theorem \ref{T2}, this yields
$$\lambda_{ij}=\int_{\bar{\cal P}_{\rm m}}
\left( \sum_{n=1}^\infty x_n{\1}_{\{i_n=j\}}\right) \nu_i({\rm d}\bx)\,,\qquad \hbox{for }i\neq j.$$
As $\sum_{j=1}^k\lambda_{ij}=0$, this entails
$$\lambda_{ii}=-\int_{\bar{\cal P}_{\rm m}}
\left( \sum_{n=1}^\infty x_n{\1}_{\{i_n\neq i\}}\right) \nu_i({\rm d}\bx)\,.$$
Thus, by \eqref{eqcons}, we obtain the general formula
\begin{equation}\label{eqintens}
\lambda_{ij}=\int_{\bar{\cal P}_{\rm m}}
\left( \sum_{n=1}^\infty x_n{\1}_{\{i_n=j\}}-\1_{\{i=j\}}\right) \nu_i({\rm d}\bx)\,,\qquad i,j\in\{1,\ldots,k\}.
\end{equation}

A slight refinement of this argument enables us to compute the finite measure
$\lambda_{ij} p_{ij}B_{ij}$. Specifically, one finds for $i\neq j$
$$\lambda_{ij}p_{ij}\int f(b)B_{ij}({\rm d} b)=\int_{\bar{\cal P}_{\rm m}}
\left( \sum_{n=1}^\infty x_n{\1}_{\{i_n=j\}}f(-\ln x_n)\right) \nu_i({\rm d}\bx)\,.$$
This gives
$$\lambda_{ij}p_{ij}(1-\hat B_{ij}(\theta))=\int_{\bar{\cal P}_{\rm m}}
\left( \sum_{n=1}^\infty {\1}_{\{i_n=j\}}(x_{n}-x_{n}^{1+\theta})\right) \nu_i({\rm d}\bx)\,.$$

Finally, the calculation of the Bernstein functions of the subordinators $S^{(i)}$ is made by reduction to the monotype situation. Specifically, we shall work under the law $\P_i$, and
 we denote by $\nu_{i}^{\dag}$ the image of 
$\nu_i$ by the map $\dag_{i}: \bar{\cal P}_{\rm m}\to {\cal P}_{\rm m}$ where $\dag_{i}(\bx)$ is the mass-partition
given by rearrangement of the terms $\1_{\{i_n=i\}}x_n$. Informally, this means that all the components of $\bx$ which are not of type $i$ are reduced to dust. Then $\nu_{i}^{\dag}$
is a (monotype) dislocation measure. It should be plain from the Poissonian construction that if we denote by $\zeta$ the instant (i.e. the first coordinate) of the first atom $(t_{i,m},\bar\pi^{(i,m)},1)$ with $\bar\pi^{(i,m)}=(\pi^{(i,m)},{\bf i})$
and $i_1\neq i$, then $\zeta$ is the first jump time of the type process $J$ and the process  killed at time $\zeta$,  $(|\Pi_1(t)|, u<\zeta)$, can be viewed as 
the process of the tagged fragment in a homogeneous  monotype fragmentation with no erosion and dislocation measure $\nu_{i}^{\dag}$. This yields 
$$\Psi^{(i)}(\theta) = \Psi^{(i,\dag)}(\theta)-\Psi^{(i,\dag)}(0)\,,$$
where, according to  Theorem 3.2 in \cite{RFC} on its page 135, 
$$\Psi^{(i,\dag)}(\theta)=\int_{{\cal P}_{\rm m}}\left(1-\sum_{n=1}^{\infty}x_n^{1+\theta}\right)\nu_{i}^{\dag}({\rm d}{\bf x})= \int_{\bar{\cal P}_{\rm m}}
\left(1- \sum_{n=1}^\infty {\1}_{\{i_n=i\}}x_{n}^{1+\theta}\right) \nu_i({\rm d}\bx)\,.$$
Hence 
$$\Psi^{(i)}(\theta) = \int_{\bar{\cal P}_{\rm m}}
\left(\sum_{n=1}^\infty {\1}_{\{i_n=i\}}(x_{n}^{1+\theta}-x_{n})\right) \nu_i({\rm d}\bx)\,.$$

Putting the pieces together in \eqref{eqK}, this establishes our claim. \QED

If we introduce
$$\underline \theta= \inf\left\{\theta\in\R: 
\int_{\bar{\cal P}_{\rm m}}\left|{\1}_{\{i=j\}}-\sum_{n=1}^\infty x_n^{1+\theta}{\1}_{\{i_n=j\}}\right|\nu_i({\rm d}{\bx})<\infty
\hbox{ for every }i,j\in\{1,\ldots,k\}\right\}\,,$$
then the  Bernstein matrix function $\theta\to {\bf \Phi}(\theta)$ 
possesses
an analytic extension to $]\underline\theta , \infty[$ 
and Theorem \ref{T3} still holds for $\theta\in ]\underline\theta , \infty[$.

\subsection{Connexion with multitype branching random walks}
The preceding analysis of the evolution of the tagged fragment provides us with
the key to shift some deep results on multitype branching random walks to homogeneous
fragmentations. The approach is quite similar to that in \cite{BeRou}, so again we shall skip details.

Just as in the monotype case, we consider the logarithms of the masses of the fragments 
and introduce for every $t\geq 0$ the 
empirical measure 
${\bf Z}^{(t)}=(Z_{1}^{(t)},\ldots, Z_{k}^{(t)})$, where
\begin{equation}\label{defZ}
Z_{j}^{(t)}\,=\,\sum_{n=1}^{\infty}\1_{\{T_{n}(t)=j\}}\delta_{-\ln
X_n(t)}\,.
\end{equation}
For every fixed step-parameter ${a}>0$, the process in discrete time
$({\bf Z}^{({a} n)}, n\in \Z_{+})$ is then a multitype branching random walk;
see \cite{BRS} for a precise definition. For the sake of simplicity, we shall focus on the case when ${a}=1$ and compute first a quantity of fundamental importance in terms of the characteristics of the fragmentation.

The analysis of multitype branching random walks relies on the Laplace transform of the intensity
$$m_{ij}(\theta)=\E_{i}\left(\int_{\R}\e^{-\theta x}Z^{(1)}_{j}({\rm d}x)\right)
=
\E_{i}\left(\sum_{n=1}^{\infty}\1_{\{T_{n}(1)=j\}} X_{n}^{\theta}(1)\right)\,.$$
Recall now that ${\bar X}(t)=|\bar\Pi(t)|^{\downarrow}$ and, from \eqref{eqsbs}, that conditionally on  ${\bar X}(t)$, the tagged fragment $|\bar\Pi_1(t)|=(\exp-S_t, J_t)$ is distributed as a size-biased sample of  ${\bar X}(t)$. Hence, for every $\theta>\underline \theta+1$, we have
$$m_{ij}(\theta)=\E_{i,0}(\exp(-(\theta-1)S_{1}), J_{1}=j)=\left(\e^{-{\bf \Phi}(\theta-1)}\right)_{ij}\,.$$

We shall now make a further assumption on the fragmentation $\bar{\bf X}$,
which will be crucial to investigate its asymptotic behavior. 
Specifically, we assume henceforth that the process $J$ of the type of the tagged fragment is ergodic, i.e. the intensity matrix ${\bf \Lambda}=(\lambda_{ij})$ 
given by \eqref{eqintens} is irreducible.
We  recall from the Perron-Frobenius theory (see for instance Seneta \cite{Sen} or Section I.6 and II.4 in Asmussen \cite{As})
that for every $\theta > \underline\theta$, the matrix $\exp(-{\bf \Phi}(\theta))$ has a unique real eigenvalue with maximal modulus which can be expressed as $\e^{-\varphi(\theta)}$.
In other words, $\varphi(\theta)$ is the eigenvalue of the Bernstein matrix ${\bf \Phi}$
with minimal real part. We also write ${\bf u}(\theta)=(u_{1}(\theta),\ldots,u_{k}(\theta))$ and 
${\bf v}(\theta)=(v_{1}(\theta),\ldots,v_{k}(\theta))$
for the left and right eigenvectors
\footnote{
It may be useful to compare our notation with that in \cite{BRS}, see in particular  Theorem 1 there.
The matrix $M(\theta)$, the eigenvalue $\rho(\theta)$ and the eigenvectors $u(\theta)$ and $v(\theta)$ there coincide respectively with  $\exp(-{\bf \Phi}(\theta-1))$,  $\exp(-\varphi(\theta-1))$, 
 ${\bf u}(\theta-1)$ and ${\bf v}(\theta-1)$ here.}
  associated with  $\e^{-\varphi(\theta)}$, normalized so that $\sum_{i=1}^{k}u_{i}(\theta)= 1$ and $\sum_{i=1}^{k}u_{i}(\theta)v_{i}(\theta)= 1$. 
  
We are now able to turn our attention to a fundamental family of martingales, which have been introduced first by Biggins in the monotype situation.

\begin{theorem}\label{T4}  Assume that the erosion coefficients ${\tt c}_{i}$ are all zero,  that the dislocation measures $\nu_{i}$ are conservative (in the sense that
\eqref{eqcons} holds), and that the intensity matrix \eqref{eqintens}
is irreducible.

\noindent {\rm (i)}The equation 
$$\varphi(\theta)=(\theta+1)\varphi'(\theta)$$
possesses a unique solution $\bar \theta \geq 0$. The function
$\theta\to \varphi(\theta)/(\theta+1)$ increases on $]\underline \theta,\bar \theta[$ and decreases on $]\bar\theta,\infty[$, and thus reaches its unique maximum at $\bar\theta$. 

\noindent {\rm (ii)} For every $\theta\in]\underline \theta,\bar \theta[$, the process
$$
M(\theta,t)= \e^{t\varphi (\theta)}\sum_{n=1}^{\infty} v_{T_{n}(t)}(\theta)X_{n}^{\theta+1}(t)\,,Ê\qquad t\geq 0$$
is a martingale which converges a.s. and in $L^1(\P_{i})$ for every type $i\in\{1,\ldots,k\}$. 
Further, this convergence is uniform for $\theta$ in any compact set in $]\underline \theta,\bar \theta[$, almost surely.
\end{theorem}

\proof (i) It can be shown that the function
\begin{equation}\label{eqeigen}
\varphi : ]\underline\theta,\infty[\to \R\hbox{ is concave, increasing, and }
\varphi(\theta)= o(\theta) \hbox{ as }\theta\to \infty.
\end{equation}
See Theorem 3.7 in Seneta \cite{Sen} for the concavity assertion, the fact that $\varphi$ increases is similar. Finally observe from Theorem \ref{T3} that $\lim_{\theta\to\infty}
\theta^{-1}{\bf \Phi}(\theta)=0$, which entails $\varphi(\theta)= o(\theta)$.
We can then follow the arguments of the proof of Lemma 1 in \cite{Be3}.

(ii) We start by recall  from the size-biased sampling formula \eqref{eqsbs} that,
if $(\f_{t})_{t\geq 0}$ denotes the natural filtration of $\bar{\bf X}$, then
\begin{eqnarray*}
M(\theta,t) &=& \e^{t\varphi (\theta)} \sum_{n=1}^{\infty}   v_{T_{n}(t)}(\theta)X_{n}^{\theta+1}(t)\\
&=& \e^{t\varphi (\theta)}\E\left( \exp(-\theta S_{t}) v_{J_{t}}(\theta)\mid \f_{t}\right)\,.
\end{eqnarray*}
As ${\bf v}(\theta)$ is a right eigenvector of $\exp(-{t\bf \Phi}(\theta))$ corresponding to the eigenvalue $\e^{-t\varphi (\theta)}$, we easily see from the Markov property of $(J,S)$  that the process $\e^{t\varphi (\theta)} \exp(-\theta S_{t}) v_{J_{t}}(\theta)$
 is a martingale in its own filtration. By projection on $(\f_{t})_{t\geq 0}$, we conclude that
 $M(\theta,t)$ is an $(\f_t)$-martingale.

By an argument of discretization analogous to that in \cite{BeRou}, 
it suffices to establish the statement when $t$ goes to infinity along the sequence
${a} n$ for some arbitrary ${a}>0$. For the sake of simplicity, we shall focus on the case
${a}=1$  without loss of generality and aim at applying Theorems 2 and 3 of \cite{BRS} to the discrete time martingale
$$
M(\theta,n)= \e^{n\varphi (\theta)}\sum_{j=1}^k v_{j}(\theta) 
\int \e^{-(\theta+1) x}Z^{(n)}_{j}({\rm d}x)
\,,Ê\qquad n\in\Z_{+}.$$

Recall from Theorem 1(ii) in \cite{BRS} that $v_{j}(\theta)\neq 0$ for every $\theta> \underline\theta$ and $j=1,\ldots,k$.
An application of the conditional Jensen's inequality to the identity
$$M(\theta,1) = \e^{\varphi (\theta)}\E\left( \exp(-\theta S_{1}) v_{J_{1}}(\theta)\mid \f_{1}\right)$$
shows that for every $\theta > \underline\theta$, there is $\alpha>1$ such that $\E_{i}(M(\theta,1)^{\alpha})<\infty$ for all types $i=1,\ldots,k$.
On the other hand,
we deduce from (i) that whenever  $\theta\in]\underline \theta,\bar \theta[$,
we can find $\alpha>1$ close to $1$ such that $\varphi(\theta)/(\theta+1)
< \varphi(\theta')/(\theta'+1)$ where $\theta'=\alpha(\theta+1)-1$.
This implies that
$$\exp(-\varphi(\alpha(\theta+1)-1)+\alpha\varphi (\theta))<1\,,$$
and Theorems 2 and 3 in \cite{BRS} now entails our claim. \QED

Let us mention an interesting consequence of Theorem \ref{T4} to the rate of decay of the largest fragment as time goes to infinity.
It follows readily from Theorem \ref{T4} that a.s.
$$\lim_{{t\to \infty}}\frac{1}{t} \ln X_{1}(t)= -\varphi'(\bar\theta)\,,$$
see e.g. Corollary 1 in \cite{Besur}. More precisely, the latter easily entails that  a.s. for every type $j\in\{1,\ldots, k\}$,
$$\lim_{{t\to \infty}}\frac{1}{t} \ln \xi_{j}(t)= -\varphi'(\bar\theta)\,,$$
where $\xi_{j}(t)=\max\{X_{n}(t), T_{n}(t)=j\}$.

\subsection{Asymptotic behavior of the empirical measure}
We shall now conclude this section by presenting a couple of applications to the asymptotic behavior of homogeneous multitype mass-fragmentations. The first belongs to the same vein as Corollary 3.3 in \cite{RFC} on its page 158. 
In the monotype case, a version of the result in discrete time goes back to Kolmogorov \cite{Ko}, in probably
the first rigorous work ever on fragmentation processes.
Roughly speaking, Kolmogorov provided an explanation to the fact which has been observed experimentally  in mineralogy, that the logarithms of the masses of mineral grains are often normally distributed.
We shall see that a similar feature holds for the more general model of multitype fragmentations.

\begin{corollary}\label{C1}  Assume that the erosion coefficients ${\tt c}_{i}$ are all zero,
 that the dislocation measures $\nu_{i}$ are conservative (in the sense that
\eqref{eqcons} holds), and that the intensity matrix \eqref{eqintens}
is irreducible. Denote for simplicity by  ${\bf u}(0)={\bf u}=(u_1,\ldots, u_k)$ the stationary distribution on $\{1,\ldots,k\}$ of the Markov chain $J$, i.e. ${\bf u}$ is the unique probability vector with
$${\bf u}{\bf \Lambda}= {\bf 0}\,.$$
Suppose  further that $\varphi$ is twice differentiable at $0$. Then
 the following limits hold in
$L^2(\P_{i})$ for any initial type $i\in\{1,\ldots, k\}$
and every  continuous bounded function
$f:\R\times\{1,\ldots,k\}\to\R$:

 $$\lim_{t\to\infty}\sum_{n=1}^{\infty}X_n(t) f(t^{-1}\ln X_n(t),T_n(t))
\,=\,\sum_{j=1}^k u_j f(-\varphi'(0),j)\,$$
and 
$$\lim_{t\to\infty}\sum_{n=1}^{\infty}X_n(t) f(t^{-1/2}(\ln X_n(t)+\varphi'(0)
t),T_n(t))\,=\,
\sum_{j=1}^k u_j\E(f({\cal N}(0,-\varphi''(0)),j))\,,$$ 
where ${\cal N}(0,-\varphi''(0))$ denotes
a centered Gaussian variable with variance $-\varphi''(0)$. 
\end{corollary}
Informally, the first limit means that the masses of most fragments decay exponentially fast
with rate $\varphi'(0)$ and their types are distributed according to the stationary law ${\bf u}$ of the Markov chain $J$. The second limit is a refinement of the first and shows that, in a pathwise sense,  fluctuations are Gaussian and independent of the type.

\proof The two limits can be established by first and second moments estimates which rely respectively on the law of large numbers and the central limit theorem for the Markov additive process $(S,J)$, and an argument of propagation of chaos. For the sake of conciseness, we shall focus on the second limit.

 Under the present assumptions,  we know from Corollary 2.8 in \cite{As} on its page 313 that as $t\to\infty$,
 $$\frac{S_t-\varphi'(0)t}{\sqrt t}\Rightarrow {\cal N}(0,-\varphi''(0))\,,$$
 where $\Rightarrow$ is used as a symbol for convergence in distribution.
On the other hand, we also have
$J(t) \Rightarrow \tau$, where $\tau$ is a random type distributed according to the stationary law ${\bf u}$. Further an easy argument using the fact that the Markov chain  $J$ mixes exponentially fast shows the asymptotic independence, in the sense that
$$\left(\frac{S_t-\varphi'(0)t}{\sqrt t}, J_t\right) \Rightarrow \left({\cal N}(0,-\varphi''(0)),
\tau\right)\,,$$
where in the right-hand side, the variables ${\cal N}(0,-\varphi''(0)$ and $\tau$ are assumed independent. 

Recall now that ${\bar X}(t)=|\bar\Pi(t)|^{\downarrow}$ and, from \eqref{eqsbs}, that conditionally on  ${\bar X}(t)$, the tagged fragment $|\bar\Pi_1(t)|=(\exp-S_t, J_t)$ is distributed as a size-biased sample of  ${\bar X}(t)$. Hence
$$\E_{i}\left(
\sum_{n=1}^{\infty}X_n(t) f(t^{-1/2}(\ln X_n(t)+\varphi'(0)
t),T_n(t))\right) = \E_{i}\left(f\left(\frac{-S_t+\varphi'(0)t}{\sqrt t}, J_t\right)\right)\,,
$$
and therefore
\begin{eqnarray*}
 & &\lim_{t\to\infty}\E_{i}\left(
\sum_{n=1}^{\infty}X_n(t) f(t^{-1/2}(\ln X_n(t)+\varphi'(0)
t),T_n(t))\right) \\
&=&\E\left(f \left({\cal N}(0,-\varphi''(0)),
\tau\right)\right)\\
&=&
\sum_{j=1}^k u_j\E(f({\cal N}(0,-\varphi''(0)),j))\,.
\end{eqnarray*}

By an argument of propagation of chaos similar to that in the proof of Corollary 3.3 in \cite{RFC} on its page 159-160, we can estimate the second moment and get
\begin{eqnarray*}
 & &\lim_{t\to\infty}\E_{i}\left(\left(
\sum_{n=1}^{\infty}X_n(t) f(t^{-1/2}(\ln X_n(t)+\varphi'(0)
t),T_n(t))\right)^2\right) \\
&=&
\left(\sum_{j=1}^k u_j\E(f({\cal N}(0,-\varphi''(0)),j))\right)^2\,.
\end{eqnarray*}
This entails the convergence in $L^2(\P_{i})$ which has been stated. \QED

Finally,  using  time discretization  techniques similar to those in \cite{BeRou}, we can translate Theorem 7 in \cite{BRS} to multitype fragmentations. This yields a pathwise large deviation limit theorem for the empirical distribution of the fragmentation which refines considerably Corollary \ref{C1}. In this direction, we shall further assume that the eigenvalue function
$\varphi$ is strictly concave and that the branching random walk ${\bf Z}^{(n)}$
is strongly non-lattice (see \cite{BRS} for the terminology), which are very natural and mild conditions.  Recall Theorem \ref{T4} and denote the terminal value of the martingale $M(\theta,t)$  by $M(\theta,\infty)$.

\begin{corollary}\label{C2} Let $h: \R\times\{1,\ldots,k\}\to\R$ be a continuous function with compact support. Under the preceding assumptions,  we have
\begin{eqnarray*}
& &\lim_{t\to\infty}  \sqrt t \, 
\e^{-t((\theta+1)\varphi'(\theta)-\varphi(\theta))}\sum_{{n=1}}^{\infty}
 h(t\varphi'(\theta)+\ln X_{n}(t),T_{n}(t))\\
&=& {M(\theta,\infty)\over \sqrt{2\pi
|\varphi''(\theta)|}}\sum_{j=1}^{k}u_{j}(\theta)\int_{-\infty}^{\infty}h(y,j)\e^{-(\theta+1)y}{\rm d}y\,,
\end{eqnarray*}
uniformly for $\theta$ in compact subsets of $]\underline \theta, \bar \theta[$, almost
surely.

 \end{corollary}

In particular, this implies that for every $a<b\in\R$,  $\theta\in]\underline \theta, \bar \theta[$
and $j\in\{1,\ldots,k\}$, there is the estimate as $t\to\infty$
\begin{eqnarray*}
& &\#\left\{n\in\N: a\e^{-t\varphi'(\theta)}\leq X_{n}(t)
\leq b\e^{-t\varphi'(\theta)} \hbox{ and } T_{n}(t)=j\right\}\\
&\sim&
  A_{\theta}
 u_{j}(\theta) t^{-1/2}\exp(t((\theta+1)\varphi'(\theta)-\varphi(\theta)))\left(\e^{-a(\theta+1)}-\e^{-b(\theta+1)}\right)\,,
\end{eqnarray*}
where $A_{\theta}$ is some strictly positive random variable with finite mean.
See Corollary 2 in \cite{BRS}.

\end{section}

\vskip 1cm
\noindent
{\bf Agradecimentos :} In August 2004, I had the pleasure to be invited in Ubatuba, to give a short course on self-similar fragmentation chains for the 8th Brazilian School of Probability.
The present paper can be viewed, in some sense, as a natural prolongation of the material which
I presented during that course. I would like to thank again Pablo Ferrari for his very kind invitation
and the wonderful organization, and Vladas Sidoravicius and Maria-Eulalia Vares
for making my stays so pleasant each time I have the chance to go to Rio.

\end{document}